\documentclass[10pt,a4paper,twocolumn]{NARMS}

\usepackage{epsfig}
\usepackage{timesmt}
\usepackage{caption}
\usepackage[backend=biber, sorting=none, style=ieee]{biblatex}
\usepackage{hyperref}
\usepackage{bbm}
\usepackage{amsfonts}
\usepackage{url}
\usepackage{array}
\usepackage{pgfplots} 
\pgfplotsset{compat=newest} 
\usepackage{multirow}
\usepackage{amssymb}
\usepackage{amsmath}
\usepackage{graphicx}
\usepackage{epstopdf}
\usepackage{color}
\usepackage[english]{babel}
\usepackage{verbatim}
\usepackage{tikz}
\usepackage{booktabs}
\usepackage{mathtools}
\usepackage{siunitx}
\usepackage{accents}

\usepackage[utf8]{inputenc}
\usepackage[english]{babel}
\usepackage{booktabs}
\usepackage{balance}
\usepackage{bbm}

\newcommand{\bm}[1]{\mbox{\boldmath $ #1 $}}

\usepackage{amsthm}

\newtheorem{remark}{Remark}
\usepackage{algorithm}
\usepackage{esvect}
\usepackage[noend]{algpseudocode}

\newtheorem{assumption}{Assumption}
\graphicspath{{figures/}}
\renewcommand{\phi}{\varphi}
\newcommand{\norm}[1]{\left\lVert#1\right\rVert}
\DeclareMathOperator*{\argmin}{arg\,min}

\confshortname{AVEC'18}
\date{July 16, 2018}
\title{Eco-driving with Learning Model Predictive Control}
\author{{Yeojun Kim, Samuel Tay, Jacopo Guanetti, Francesco Borrelli} \\
{\aff{Department of Mechanical Engineering}} \\
{\aff{University of California, Berkeley, USA}} \\
{\aff{\normalfont\normalsize E-mail: yk4938@berkeley.edu}}\\
{\aff{\normalfont\normalsize Topics/model predictive control, learning, fuel efficiency}}
}

\abstract{
We present a predictive cruise controller which iteratively improves the fuel economy of a vehicle traveling along the same route every day. Our approach uses historical data from previous trip iterations to improve vehicle performance while  guaranteeing a desired arrival time. The proposed predictive cruise controller is based on the recently developed Learning Model Predictive Control (LMPC) framework, which is extended in this paper to include time constraints. Moreover, we reformulate the modified LMPC with time constraint into a computationally tractable form. Our method is presented in detail, applied to the predictive cruise control problem, and validated through simulations.
}

\begin{document}
	
\maketitle

\section{Introduction}\label{sec:Introduction}
The U.S. Department of Transportation estimates that nearly 120 million Americans commute each day with an average commute time of 26 minutes\cite{cohen2007making}. In addition, 90 $\%$ of Americans drive to work with privately owned vehicles \cite{dews2013ninety}. These daily commutes often follow the same route.

Advanced driving assistant systems (ADAS) aim to assist drivers in performing common driving tasks and maneuvers safely \cite{knoll2014driving, carvalho2015automated}. With recent advancements in perception technologies and computational power, ADAS can play an important role to improve the energy efficiency of a vehicle, as well as its safety \cite{turri2017model, guanetti2018control}.  
Predictive cruise control (PCC) is an example of this potential \cite{van2009predictive,kohut2009integrating}.
PCC is a longitudinal velocity control driving assistance system that uses look-ahead information about the downstream road.
This information includes static information such as speed limits or road grade, and dynamic information such as traffic speed or intersection delays.

Often PCC is cast as an optimization problem. 
In \cite{lattemann2004predictive, turri2017cooperative}, optimization methods to calculate the optimal velocity trajectory for maximizing energy efficiency are presented.
In \cite{asadi2011predictive, sun2018robust}, an optimal control problem is formulated to minimize braking when the vehicle is going through a series of traffic lights.
However, their approaches are limiting in the following senses. 
First, they require a priori knowledge about the environment. Second, the complexity of their approaches increases with the length of trip. Third, they cannot ensure global safety and time constraint guarantees unless the horizon of the optimal control problem is long enough.

In this work we tackle the problem of PCC using the Learning Model Predictive Control (LMPC) framework presented in \cite{rosolia2017learning}. 
This reference-free controller is an attractive approach because it is able to solve long horizon optimal control problems while ensuring global safety constraint guarantees. 
It has been successfully implemented in autonomous driving applications \cite{rosolia2017autonomous, brunner2017repetitive}. 
However, the original work of LMPC does not address the problem of completing the task within a given time limit as it is usually formulated for minimum time problems, where this issue does not arise. As a result, depending on how the cost of the optimal control problem is designed, the task may result in taking more time  to actually finish a task.

In order to apply the LMPC framework to PCC of a vehicle which repeatedly drives along the same route, we modify the original LMPC so that we can enforce the task to be completed within a time limit.
Moreover, in every trip, our controller can learn the static environment features such as road grade and attempts to improve its performance, which includes total fuel consumption and comfort.
The contributions of this paper can be summarized as follows:
	\begin{itemize}
	 \item We present the new design for LMPC with time constraint which accommodates a constraint on the total duration of the task.
	 \item We reformulate LMPC with time constraint into a computationally tractable formation so that it can be implemented in real time.
	 \item We design a predictive cruise control using the proposed LMPC with time constraint and show its effectiveness with simulation results.
	\end{itemize}

The remainder of this paper is organized as follows. Section 2 introduces the vehicle dynamics model and a fuel estimation model, and defines the optimization problem for the predictive cruise control. Section 3 formulates the learning model predictive controller with time constraint. Section 4 provides a tractable reformulation of LMPC with time constraint for real time implementation. Section 5 demonstrates the effectiveness of our controller design as a predictive cruise controller for a repetitive trip in the Berkeley area. Section 6 concludes the work with some remarks.

\section{Problem Formulation} \label{sec:problem_formulation}
We aim to build the predictive cruise control of a vehicle performing repetitive trips along a fixed route, subject to a position-dependent road grade and a total completion time constraint. We refer to each successful trip of the route as an \textit{iteration}. Each iteration has the same boundary conditions (initial and final position, speed, and force), and has to be completed within the desired time limit $N_f$; $x_0=x_s$ and $x_{N_f}=x_f$. At each iteration, the controller finds a velocity trajectory that maintains or improves the specified vehicle performance objective (such as fuel economy or comfort).

\subsection{Vehicle Dynamics and Fuel Estimation} \label{sec:vehicle_dynamics}
The simplified vehicle longitudinal dynamics are modeled as a connection of the first order systems with parameters identified from experiments. The states of our model include the distance travelled along the route $s$, the velocity $v$, and the force at the wheel $F$. The inputs to the model are $F_{t}$ and $F_{b}$, which represent the wheel-level desired traction and braking forces, respectively; $u_k = \begin{bmatrix} F_{t,k} & F_{b,k} \end{bmatrix}^T$. Denoting the state at time $k$ as $x_k = \begin{bmatrix} s_{k} & v_{k} & F_{k} \end{bmatrix}^T$, the system dynamics are
\begin{equation}\label{vehicle_model}
x_{k+1} = \left[ {\begin{array}{cc}
   s_k + t_s v_k\\
   v_k + \frac{t_s}{m} (F_k + F_{R,k}) \\
  \left(1- \frac{t_s}{\tau}\right)F_k + \frac{t_s}{\tau}(F_{t,k} + F_{b,k})\\
  \end{array} } \right]
\end{equation}
where $t_s$ is the sampling time; $m$ and $\tau$ are the mass and the time constant of force actuation, respectively; $F_{R}$ is the resistance force which includes aerodynamic drag, rolling resistance, and gravitational force. This force can be represented as
\begin{align} \label{eq:resistance_force}
F_{R,k} = m g c_r \textnormal{cos}(\theta) + m g \textnormal{sin}(\theta) + \frac{1}{2} \rho A C_d v_k^2\,
\end{align}
where $\theta$ is the road pitch angle; $g$ and $c_r$ are the gravity and rolling coefficients, respectively; $\rho$, $A$, and $C_d$ are the air density, frontal area, and drag coefficient, respectively. 
Moreover, the road angle is approximated using a quadratic function of the distance 
\begin{align}\label{eq:angle_eqn}
\hat{\theta}_t(s_k) = a_{0,k} + a_{1,k} s_k + a_{2,k} s_k^2\, ,
\end{align}
where $(a_{0,k}, a_{1,k}, a_{2,k})$ are parameters that are computed as follows.
During each iteration of the trip, we store velocity and force values at each position along the route; we then introduce $\bar{\theta}$ to estimate the angle by inverting the dynamics \eqref{vehicle_model}. At each time step $k$ of the $j$-th iteration, given $s_k^j$ the parameters $(a^j_{0,k}, a^j_{1,k}, a^j_{2,k})$ are estimated on-line solving the following least mean squares problem,
\begin{equation} \label{eq:angle_app}
\displaystyle \argmin_{(a^j_{0,k}, a^j_{1,k}, a^j_{2,k})} \sum_{(p,l) \in \mathcal{G}(s_k^j)}  \norm{\begin{bmatrix}  1 & s^l_{p} & (s^l_{p})^2\end{bmatrix} \begin{bmatrix}  a^j_{0,k} \\a^j_{1,k} \\ a^j_{2,k}\end{bmatrix}  - \bar{\theta}^l_p}\\
\end{equation}
where $\mathcal{G}(s_k^j)$ is the set of indices and the iteration numbers with the following property
\begin{equation}
\mathcal{G}(s_k^j) = \Big\{(p,l): s^{j}_k \leq s^{l}_p \leq s^{j}_k + d_{\textnormal{f}}\Big\}.
\end{equation}
where $d_f$ is look-ahead distance which is considered a tuning  parameter. A similar method is adopted for system identification of road curvature in \cite{rosolia2017autonomous}.

In the remainder of this paper, the approximated vehicle dynamics model \eqref{vehicle_model}-\eqref{eq:angle_app} is compactly rewritten as 
\begin{equation}
    x_{k+1} = f(x_k, u_k).
\end{equation}
Also, we consider state and input constraints in the form
	\begin{subequations}\label{stateinputconstraints}
	\begin{align}
	x_k \in \mathbb{X}&:= \{(s,v,F): 0\leq v \leq v_{\textnormal{max}}, \, \nonumber \\ & \quad \quad \quad \quad \quad F_{\textnormal{min}} \leq F \leq F_{\textnormal{max}}  \},\\
	u_k \in \mathbb{U}&:= \{(F_{t}, F_{b}): 0\leq F_{t} \leq F_{\textnormal{max}},\, F_{\textnormal{min}} \leq F_{b} \leq 0 \}.
	\end{align}
	\end{subequations}

We are interested in designing a predictive cruise controller which tries to improve its performance as it repeats the same route. The performance objective can be fuel consumption, jerk, or travel time. In this paper we seek a better fuel economy by improving both so-called \textit{tank-to-wheels} and \textit{wheels-to-miles} efficiency \cite{guzzella2007vehicle}. Higher \textit{tank-to-wheels} and \textit{wheels-to-miles} efficiency involve improving peak management of engine, shaping velocity, and reducing the aerodynamic and rolling losses. In order to approximate the fuel consumption, we adopted a polynomial fuel approximation method from \cite{kamal2010board} which can be written the following form,
\begin{align}\label{eq:fuel_eqn}%
f_{\textnormal{fuel}}(v, F_t) = f_{\textnormal{cruise}}(v) + f_{\textnormal{accel}}(v, F_t)\, ,
\end{align}%
where
\begin{align}
&f_{\textnormal{cruise}}(v) = (b_0 v + b_1 v^2 + b_2 v^3 )\, , \nonumber\\
&f_{\textnormal{accel}}(v, F_t) = F_t(c_0 + c_1 v + c_2 v^2)\, , \nonumber
\end{align}%
and $(b_0,b_1,b_2, c_0,c_1,c_2)$ are parameters identified by least mean squares fitting of the experimental fuel rate data. The goal of our controller is to minimize this estimate of the fuel consumption, $h(x,u) = f_{\textnormal{fuel}}(v, F_t)$. 

\subsection{Predictive Cruise Control Problem}
\label{sec:PCC_problem}
For each iteration of trip, we can formulate the predictive cruise control problem as the following constrained finite horizon optimal control problem.
\begin{subequations}\label{pcc_problem_form}
\begin{align}
    &\displaystyle \min_{u^j_{0},...,u^j_{N_f-1}}  \sum_{k=0}^{N_f-1} h(x^j_{k},u^j_{k}) \label{stage_cost_pcc}\\
	& \textnormal{subject to} \nonumber \\
	&\;\;\, x^j_{0} = x_s, \;\;\, x^j_{N_f} = x_f \label{init_condition_pcc}\\
	&\;\;\, x^j_{k+1} = f(x^j_{k}, u^j_{k}), \quad \forall k \in [0,...,N_f-1], \label{system_update_pcc} \\
	&\;\;\, x^j_{k} \in \mathbb{X}, \, u^j_{k} \in \mathbb{U}, \quad \forall k \in [0,...,N_f-1], \label{system_cons_pcc}
\end{align}
\end{subequations}%
where $j$ is the iteration number; \eqref{init_condition_pcc} and \eqref{system_update_pcc} represent the boundary conditions and the vehicle dynamics, respectively; \eqref{system_cons_pcc} represents the state and input constraints; The stage cost, $h(\cdot)$ in \eqref{stage_cost_pcc}, represents the estimated fuel consumption.

%% end space

\section{Learning Model Predictive Control with Time Constraint}\label{sec:lmpc_tc}

In this section, a formulation of LMPC with time constraint is proposed. Solving a finite time constrained optimal control problem such as \eqref{pcc_problem_form} in real time can be difficult, especially when $N_f$ is large. 
Therefore, we design LMPC which tries to solve the problem \eqref{pcc_problem_form} and can be implemented in real time. 
In previous works, LMPC was introduced for repetitive and iterative tasks \cite{rosolia2017learning}. 
LMPC leverages past data to progressively improve performance while ensuring recursive feasibility, asymptotic stability, and non-increasing cost at every iteration. 
In this work, we extend the LMPC framework with a constraint on the time required to complete the task. 
In other words, we guarantee that each iteration or repetition does not exceed a total time limit.

\begin{remark} \label{rm1}
	 In the original work of LMPC in \cite{rosolia2017learning}, the optimal control problem is defined on infinite horizon. In our problem, we focus on the finite time formulation \eqref{pcc_problem_form}.
\end{remark}

\subsection{Time Sampled Safety Set} \label{sec:SS_timed}
We denote the input sequence applied to the dynamics \eqref{vehicle_model} and the corresponding closed loop state trajectory at j-th iteration as
\begin{subequations}   
\begin{align}
  \bm{u^{j}} & = [u^{j}_{0}, u^{j}_{1}, ...,  u^{j}_{N_f-1}], \label{input_traj}\\
  \bm{x^{j}} & = [x^{j}_{0}, x^{j}_{1}, ...,  x^{j}_{N_f}]  \label{state_traj}
\end{align}
\end{subequations}
where $u^j_t$ and $x^j_t$ are the input and the state at time $t$ of the j-th iteration, respectively.

The main contribution of LMPC with time constraint is the modification of the safety set in \cite{rosolia2017learning} to the \textit{Time Sampled Safety Set}. We define the time sampled safety set $\mathcal{SS}_{\textnormal{time}}^j$ at j-th iteration as
\begin{align} \label{eq:safety_set}
\mathcal{SS}_{\textnormal{time}}^{j}(t) = \displaystyle \cup_{i=1}^{j} \displaystyle \cup_{k=t}^{N_f} \, x_k^i 
\end{align}
where $N_f$ is the time limit for each iteration. 
The difference from the original definition of the safety set is that it only includes the states visited during the remaining time, $N_f - t$.
Note that $\mathcal{SS}_{\textnormal{time}}^{j}(N_f)$ is only $x_f$ since each trip must finish within the time constraint $N_f$.

\subsection{Preliminaries} \label{sec:Preliminaries}
In this section, we introduce some terminology used for the LMPC problem with time constraint.

At time $t$ of the j-th iteration, we define the cost-to-go associated with the input sequence \eqref{input_traj} and the corresponding state trajectory \eqref{state_traj} as
\begin{equation}  \label{eq:Jcost}
    J^{j}_{t \rightarrow N_f}(x_{t},t) = \sum_{k=t}^{N_f} h(x_{k}^{j},u_{k}^{j})
\end{equation}
where $h(\cdot)$ is the stage cost function such as the fuel estimation function \eqref{eq:fuel_eqn}. We have the following assumption about the stage cost $h(x,u)$.
\begin{assumption} \label{as:1}
$h(\cdot, \cdot)$ is a continuous function which has the following property:
\begin{align} 
h(x_f, u) = 0 \textnormal{ and } h(x, u)\geq 0 \,\, \forall x \in \mathbb{R}^{n_x}, \,  u \in \mathbb{R}^{n_u}. \nonumber
\end{align}
where $n_x$ and $n_u$ are the dimensions of $x$ and $u$, respectively.
\end{assumption}

For any $x\in \mathcal{SS}^j_{\textnormal{time}}(t)$, we can define the minimum cost-to-go function $Q^{j}(x,t)$ as
\begin{equation}
\begin{aligned}
    Q^{j}(x,t)& = \\
    &\begin{cases}
       \displaystyle \min_{(i,l)\in\mathcal{F}^j(x,t)}  J_{l\rightarrow N_f}^i(x, l) \text{ if }x\in\mathcal{SS}^j_{\textnormal{time}}(t)  \\
       \, + \infty \quad \textnormal{else}
       \end{cases}
\end{aligned}
\end{equation}
where $\mathcal{F}^j(x,t)$ is defined as
\begin{subequations}
    \begin{align}
    \mathcal{F}^j(x,t) = \{ (i,l): &i\in [0,j], \, l \geq t, \, x= x^i_l \nonumber \\
    &\textnormal{ such that }  x^i_l \in\mathcal{SS}_{\textnormal{time}}^{i}(l) \}.   \nonumber \end{align}
\end{subequations}
Note that the definition of the function $Q(\cdot,\cdot)$ is modified from the original definition in \cite{rosolia2017learning} because we use the new time sampled safety set $\mathcal{SS}_{\textnormal{time}}^{i}(\cdot)$.

\subsection{LMPC with Time Constraint Formulation} \label{sec:prob_form}
At time $t$ of iteration $j \geq 1$, our LMPC with time constraint solves the following optimization problem:
\begin{subequations}\label{lmpc_problem}
\begin{align}
    &\displaystyle \min_{u_{0|t},...,u_{N-1|t}}  \sum_{k=0}^{N-1} h(x^j_{k|t},u^j_{k|t}) + Q^{j-1}(x_{N|t},t) \\
	& \textnormal{subject to} \nonumber \\
	&\;\;\, x^j_{0|t} = x^j_t,   \label{init_condition}\\
	&\;\;\, x^j_{k+1|t} = f(x^j_{k|t}, u^j_{k|t}), \quad \forall k \in [0,...,N-1], \label{system_update} \\
	&\;\;\, x^j_{k|t} \in \mathbb{X}, \, u^j_{k|t} \in \mathbb{U}, \quad \forall k \in [0,...,N-1], \label{system_cons}\\
	&\;\;\, \mathbbm{1}_A (k) x^j_{k|t} = \mathbbm{1}_A (k) x_f, \quad \forall k \in [0,...,N-1], \label{stay_cons} \\
	&\;\;\, x^j_{N|t} \in \mathcal{SS}_{\textnormal{time}}^{j-1}(t+N) , \label{term_cons}
% 	&\;\;\, x^i_{N|t} \in \mathcal{SS}_{\textnormal{time}}^{i-1}(t+N - \mathbbm{1}_A(t+N)* (t+N-N_f)) \label{term_cons}
\end{align}
\end{subequations}%
where $\mathbbm{1}_A(\cdot)$ is the indicator function of the set $A$ defined as $A = \{t\in\mathbb{R}:t\geq N_f\}$.
Constraints \eqref{init_condition} and \eqref{system_update} represent the initial condition and vehicle dynamics, respectively; \eqref{system_cons} represents the state and input constraints; \eqref{stay_cons} is the constraint which forces the system to stay at $x_f$ at $t\geq N_f$; \eqref{term_cons} is the terminal constraint which imposes the system to be driven into the safe set sampled from last iteration.
$Q^{0}(\cdot,\cdot)$ and $\mathcal{SS}_{\textnormal{time}}^{0}(\cdot)$ are defined by the initial successful trip.

The resulting optimal states and inputs of \eqref{lmpc_problem} are denoted as
\begin{subequations}   
\begin{align}
  \bm{x^{*,i}_t} & = [x^{*,i}_{0|t}, x^{*,i}_{1|t}, ...,  x^{*,i}_{N|t}],\\
  \bm{u^{*,i}_t} & = [u^{*,i}_{0|t}, u^{*,i}_{1|t}, ...,  u^{*,i}_{N-1|t}] \label{optimal_input}.
\end{align}
\end{subequations}
Then, the first input $u^{*,i}_{0|t}$ is applied to the system during the time interval $[t,t+1)$; 
\begin{equation}  \label{lmpc_input} 
      u^i_t = u^{*,i}_{0|t}.
\end{equation}
At the next time step $t+1$, a  new optimal  control  problem  in  the  form of \eqref{lmpc_problem}, based on new measurements of the state, is solved over a shifted horizon, yielding a \textit{moving} or \textit{receding} horizon control strategy with control law.

% \begin{remark} \label{rm1}
% 	 LMPC with time constraint \eqref{lmpc_problem}-\eqref{lmpc_input} can be modified to include dynamic constraints with minor modifications.
% \end{remark}

It is noted that with the assumptions \ref{as:1}, LMPC with time constraint \eqref{lmpc_problem}-\eqref{lmpc_input} is recursively feasible and the cost of each iteration monotonically decreases. The proof is similar to the original work of LMPC in \cite{rosolia2017learning}. The key difference between the two LMPC frameworks is that the time sampled safety set shrinks in the course of time whereas in \cite{rosolia2017learning}, the safety set is time independent; however, this doesn't affect the proof because the time sampled safety set always includes at least one point which guarantees the existence of the feasible input and the cost decrease at the next time step; therefore, iteration cost decreases as the iteration progresses. 

\begin{remark} \label{rm2}
	 LMPC with time constraint \eqref{lmpc_problem}-\eqref{lmpc_input} can be reformulated as a robust control and take into account of dynamic environment with minor modifications \cite{rosolia2017robust}.
\end{remark}

% The proof for Remark \ref{rm1} is

% \subsection{Properties of LMPC with Time Constraint} \label{sec:recursv_Feas_cost_decreas}
% Following the assumptions from \cite{rosolia2017learning} that $x_f \in $ is an equilibrium point and there is no model mismatch, we can guarantee our formulation can guarantee... 
 
%  keep it short and refer to ugo paper for detailed? or maybe delete.

%% space between sections

%% end space

\section{LMPC Relaxation for Predictive Cruise Control}\label{sec:lmpc_relaxation}
In this section we apply the LMPC with time constraint \eqref{lmpc_problem}-\eqref{lmpc_input} to a predictive cruise controller subject to repetitions of the same commute with a total time limit. Because solving the optimization problem \eqref{lmpc_problem} in real time is computationally challenging \cite{rosolia2017autonomous}, we use an approximation method for \eqref{lmpc_problem}: we introduce approximation functions for the time sampled safety set, $\mathcal{SS}_{\textnormal{time}}(\cdot)$, and the terminal set, $Q(\cdot)$.

Because we restrict the vehicle velocity to be positive semi-definite, the distance travelled $s$ always monotonically increases with time $t$. Therefore, we use $s$ to shrink the time sampled set at each time. 
At time $t$ of j-th iteration, we approximate $\mathcal{SS}_{\textnormal{time}}(\cdot)$ with 
\begin{equation} \label{eq:approx_safety_set}
\begin{aligned}
\hat{\mathcal{SS}}_{\textnormal{time}}(t) = \Big\{ (s,v,F) : &s \geq s^{j-1}_t, \\
& \begin{bmatrix} v \\ F\end{bmatrix}  =  \Lambda \begin{bmatrix} 1 \\ s \\ s^2\end{bmatrix} \Big\}
\end{aligned}
\end{equation}
where
$\Lambda \in \mathbb{R}^{2\times3}$ is the solution of the following least mean square optimization problem
\begin{equation}
\displaystyle \argmin_{\Lambda}  \sum_{k \in \mathcal{T}(s_t^j)} \norm{\begin{bmatrix} v^{j-1}_k \\ F^{j-1}_k \end{bmatrix} - \Lambda \begin{bmatrix} 1 \\ s^{j-1}_k \\ (s^{j-1}_k)^2\end{bmatrix}}\\
\end{equation}
where $\mathcal{T}(s_t^j)$ defines the time steps in which the distance travelled during the previous $j-1$-th iteration of trip is between the current distance and a far enough distance forward, $d_f$:
\begin{equation}
\mathcal{T}(s_t^j) = \Big\{k: s^{j}_t \leq s^{j-1}_k \leq s^{j}_t + d_{\textnormal{f}} \Big\}.
\end{equation}
It's noted that $\hat{d}_f$ is a tuning parameter decided by the control designer.

In order to approximate the cost-to-go function $Q(\cdot, \cdot)$, we introduce the third-order polynomial function $\mathcal{C}(\cdot)$
\begin{equation}
\mathcal{C}(s) = \begin{bmatrix} 1 & s & s^2 & s^3 \end{bmatrix} \Delta 
\end{equation}
where
$\Delta \in \mathbb{R}^4$ is the solution of the following least mean square optimization problem
\begin{equation}
\displaystyle \argmin_{\Delta}  \sum_{k \in \mathcal{T}(s_t^j)} \norm{J^{j-1}_{k\to N_f}-\begin{bmatrix} 1 & s_k^j & (s_k^j)^2 & (s_k^j)^3 \end{bmatrix} \Delta }
\end{equation}
where $J^{j-1}_{k\to N_f}$ is the cost-to-go function defined in \eqref{eq:Jcost}.

Finally, we approximate $Q(\cdot, \cdot)$ with
\begin{equation} \label{eq:Q_app}
\begin{aligned}
    \hat{Q}^{j}(x,t) =& \\
    &\begin{cases}
       \mathcal{C}(s) \text{ if }x\in \hat{\mathcal{SS}}^j_{\textnormal{time}}(t)  \\
       \, + \infty \quad \textnormal{else}
       \end{cases}
\end{aligned}
\end{equation}
where $x = \begin{bmatrix} s & v & F \end{bmatrix}^T$.

\subsection{LMPC with Time Constraint Relaxation} \label{sec:lmpc_pcc}
With the approximation functions \eqref{eq:approx_safety_set}-\eqref{eq:Q_app}, we can reformulate LMPC with time constraint \eqref{lmpc_problem}-\eqref{lmpc_input} as the following optimal control problem:
\begin{subequations}\label{lmpc_problem_approx}
\begin{align}
    &\displaystyle \min_{u_{0|t},...,u_{N-1|t}}  \sum_{k=0}^{N-1} h(x^j_{k|t},u^j_{k|t}) + \hat{Q}^{j-1}(x_{N|t},t) \\
	& \textnormal{subject to} \nonumber \\
	&\;\;\, x^j_{0|t} = x^j_t,  \\
	&\;\;\, x^j_{k+1|t} = f(x^j_{k|t}, u^j_{k|t}), \quad \forall k \in [0,...,N-1],  \\
	&\;\;\, x^j_{k|t} \in \mathbb{X}, \, u^j_{k|t} \in \mathbb{U}, \quad \forall k \in [0,...,N-1], \\
	&\;\;\, \mathbbm{1}_A (k) x^j_{k|t} = \mathbbm{1}_A (k) x_f, \quad \forall k \in [0,...,N-1], \\
	&\;\;\, x^j_{N|t} \in \hat{\mathcal{SS}}_{\textnormal{time}}^{j-1}(t+N) , 
% 	&\;\;\, x^i_{N|t} \in \mathcal{SS}_{\textnormal{time}}^{i-1}(t+N - \mathbbm{1}_A(t+N)* (t+N-N_f)) \label{term_cons}
\end{align}
\end{subequations}%
The resulting optimal states and inputs of \eqref{lmpc_problem_approx} are denoted as
\begin{subequations}   
\begin{align}
  \bm{\hat{x}^{*,i}_t} & = [x^{*,i}_{0|t}, x^{*,i}_{1|t}, ...,  x^{*,i}_{N|t}],\\
  \bm{\hat{u}^{*,i}_t} & = [u^{*,i}_{0|t}, u^{*,i}_{1|t}, ...,  u^{*,i}_{N-1|t}] \label{optimal_input_approxi}.
\end{align}
\end{subequations}
Then, the first input $u^{*,i}_{0|t}$ is applied to the system during the time interval $[t,t+1)$; 
\begin{equation}\label{eq:lmpc_input_approx} 
      \hat{u}^i_t = u^{*,i}_{0|t}.
\end{equation}
At the next time step $t+1$, a  new optimal  control  problem \eqref{lmpc_problem_approx} with new measurements of the state, is solved over a shifted horizon.

%% space between sections

%% end space

\section{Simulation Results}\label{sec:sim}
\begin{figure}[!h]
	\centering
 	\includegraphics[width=\linewidth,height=\textheight,keepaspectratio]{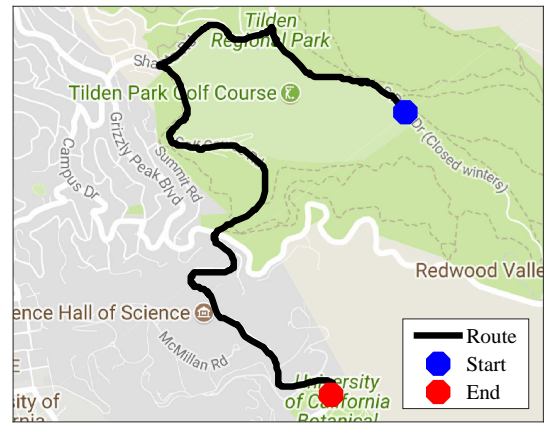} \vspace{-20pt}
	\caption{Fixed route from an origin A and a destination B in the Berkeley hill area}
	\label{fig:map}
\end{figure}

\begin{figure}[!h]
	\centering
 	\includegraphics[width=\linewidth,height=\textheight,keepaspectratio]{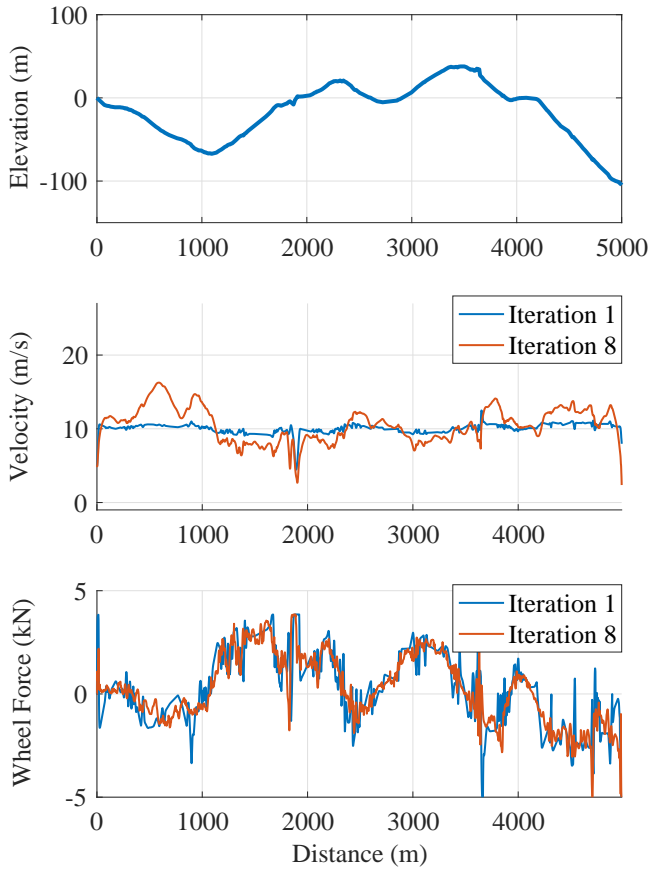}
	\caption{Plots of slope (top), closed loop velocity trajectories (middle), and closed loop wheel force trajectories (bottom) }
	\label{fig:trajectory}
\end{figure}

In this section we validate the proposed LMPC controller \eqref{lmpc_problem_approx}-\eqref{eq:lmpc_input_approx} with simulation results. A vehicle is repeating the same trip from $s_s = 0 \si{\m}$ to $s_f = 5000\si{\m}$ in the Berkeley hills area, depicted in Figure~\ref{fig:map}. This route is subject to position-dependent slope as shown in the top plot in Figure~\ref{fig:trajectory}. Each trip is initialized with $x_s = [s_s, 0, 0]$ and ends with $x_f = [s_f, 0, 0]$. We initialize the trip with a simple velocity tracking controller with a constant velocity reference. 

Figure \ref{fig:trajectory} depicts the closed loop trajectories for the first iteration of the trip and the 8-th iteration of the same trip. In every iteration of the trip, the arrival times does not exceed the terminal time. Also, as the iteration progresses, the velocity becomes higher in downhill sections and lower in uphill sections. This trend helps decrease the total fuel consumption of the trip as it uses the downhill regions to speed up and the uphill regions to slow down, leading to reduced acceleration.   
This behavior is also seen in force trajectories. Over the course of iterations, only in uphills, our controller maintains positive wheel force whereas in downhills, it tends to apply less braking (except near the end of trip where the vehicle must come to a full stop); therefore, it wastes less amount of energy.

\begin{figure}[!h]
	\centering
 	\includegraphics[width=\linewidth,height=\textheight,keepaspectratio]{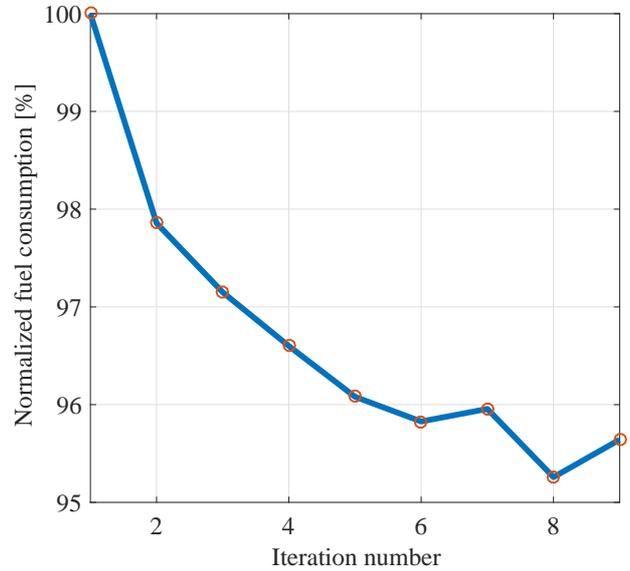} \vspace{-20pt}
	\caption{Plot of fuel consumption per iteration}
	\label{fig:cost_plot}
\end{figure}

Figure \ref{fig:cost_plot} depicts the normalized 
total cost (fuel consumption) for each iteration of the complete trips when the first trip is completed with a constant velocity tracking controller. As seen, the total cost generally decreases as the iteration number increases. There is about $4.5\%$ reduction in fuel consumption only after 8-th trip compared to the 1st trip.
It is also noted that the learning rate decreases with the iterations, as the total cost converges.
This result is analogous to those in other applications of LMPC \cite{rosolia2017autonomous,brunner2017repetitive}.

%% space between sections

%% end space 

\section{Conclusion}\label{sec:conclusion}
In this paper we proposed a predictive cruise controller which improves fuel efficiency by learning from the historical trips. The key aspect is the modification of learning model predictive cruise control to guarantee completion of the task within a total time constraint, while still improving the control performance. We validated our controller with simulation results. Future work includes the experimental validation of our developed predictive cruise controller and its modification to account for disturbances from the front vehicle. Finally, we can apply the learning model predictive control with time constraint to other real world problems.  

\section{Acknowledgement}
The information, data, or work presented herein was funded in part by the Advanced Research Projects Agency-Energy (ARPA-E), U.S. Department of Energy, under Award Number DE-AR0000791. The views and opinions of authors expressed herein do not necessarily state or reflect those of the United States Government or any agency thereof.

%% space between sections

%% end space

\printbibliography
\end{document}